\documentclass[reqno]{amsart}

\usepackage{amssymb}
\usepackage[initials]{amsrefs}
\usepackage[all]{xy}
\usepackage{mathrsfs} 
\usepackage{graphicx} 

\newtheorem{theorem}{Theorem}[section]
\newtheorem{proposition}[theorem]{Proposition}
\newtheorem{lemma}[theorem]{Lemma}
\newtheorem{corollary}[theorem]{Corollary}
\newtheorem{definition}[theorem]{Definition}
\newtheorem{assumption}[theorem]{Assumption}
\newtheorem{remark}{Remark}

\def\ainf{A_{\infty}}
\def\tens{\otimes}   
\def\Hom{\text{Hom}}
\def\mA{\mathcal{A}} 
\def\mB{\mathcal{B}} 
\def\sF{\mathscr{F}} 
\def\sG{\mathscr{G}} 
\def\ga{\gamma}       
\def\Z{\mathbf{Z}}   
\def\R{\mathbf{R}}   
\def\g{\mathfrak{g}} 
\def\h{\mathfrak{h}} 
\def\j{\mathfrak{w}} 
\def\a{\alpha}       
\def\id{\text{Id}}   
\def\at{\alpha_t}    
\def\ha{\hat{\alpha}} 
\def\hG{\hat{\mathscr{G}}} 
\def\E{\underline{E} } 
\def\lR{\underline{\textbf{R}}} 
\def\der#1#2{\frac{\partial #1}{\partial #2}} 
\def\br#1#2{\left[#1,#2\right]} 
 
\begin{document}

\title{Families of $\ainf$ algebras and homotopy group actions}

\author{Emma Smith Zbarsky}
\email{emmas@math.uchicago.edu}
\address{5734 S. University Avenue
         Chicago, Illinois 60637
         USA}
\thanks{The author was supported in part by a Bell Labs Graduate Research Fellowship}
      



\begin{abstract}
We define homotopy group actions in terms of families
of $\ainf$ algebras indexed by a manifold $M$. We give explicit
formulae for the $\ainf$ morphism induced by a path on the manifold
and for the $\ainf$ homotopy corresponding to a pair of homotopic
paths. Finally, we compute examples for finite groups and finitely
generated free nonabelian groups and determine that every homotopy
group action by a finite group is homotopic to a strict group action.
\end{abstract}


\maketitle


\section{Introduction}

As a postdoctoral fellow, James Stasheff studied group-like
topological spaces in \cite{hspaces} building on work by Sugawara in
\cite{sugawara}. He began by defining the concept of an $\ainf$
space. Initially these ideas were found useful in homotopy
theory. Generalizations were constructed including Boardman and Vogt's
machinery of topological PROPs \cite{boardvogt}, May's introduction of
operads \cite{may} and Adams' discussion of infinite loop spaces
\cite{adams}. In the nineties, $\ainf$ structures were found to have
significant presence in deformation theory, topology, and physics,
with \cite{gj}, \cite{fukaya93}, \cite{stasheff}, and
\cite{penkavaschwarz}, while Stasheff's birthday conference
contributed \cite{mccleary}. Building off of Fukaya's work, Kontsevich
conjectured homological mirror symmetry in a talk at the 1994 ICM
\cite{kontsevich}. Several special cases of homological mirror
symmetry have since been proven, notably by Polishchuk and Zaslow in
\cite{polzaslow}, Seidel in \cite{hmsqs,hmsg2c} and Abouzaid and Smith
in \cite{hms4t}. Partial proofs in other cases have been given by
Kontsevich and Soibelman in \cite{kontsoibel}, and Fukaya in
\cite{fuk02}. In a Fukaya-Seidel category, because the $\ainf$
structure arises from intersecting Lagrangians on a symplectic
manifold it is natural to wonder how group actions on the manifold may
affect the $\ainf$ structure.

\subsection{Definitions, conventions, and notation}

 An \emph{$\ainf$ algebra $\mA$} consists of a graded $K$-module $V$
 together with a sequence of maps $\mu_{\mA}^k: V^{\tens k} \to
 V[2-k]$, $k \geq 0$ that satisfy the sequence of relations
  \begin{equation}\label{ainfstruc}\sum_{k+r-1=n}\sum_{j=0}^{k-1}
    (-1)^{\maltese_j}\mu^k_{\mA}(a_1\tens \cdots \a_j\tens
    \mu^r_{\mA}(a_{j+1}\tens \cdots \tens a_{j+r}) \tens
    a_{j+r+1}\tens \cdots \tens a_n)=0, \end{equation} for $ n \geq
    0,$ where $\maltese_j = \sum_{i=1}^j (|a_i|-1)$. Now
  let $\mB$ be an $\ainf$ algebra with underlying $K$-module $W$. An
  \emph{$\ainf$ morphism $\sF: \mA \to \mB$} consists of a sequence of
  maps \[\sF^n: V^{\tens n} \to W[1-n], \hspace{.2in} n \geq 1\] which
  satisfy the corresponding sequence of relations
\begin{align}\label{ainfmorphrel}
&\sum_{n \geq r\geq 1 \atop{s_1 + \cdots s_r = n}}
  \mu^{r}_{\mB}(\sF^{s_1}(a_1,\ldots, a_{s_1}),
  \ldots,\sF^{s_r}(a_{n-s_r+1}, \ldots, a_n)) = \\ \nonumber
  &\hspace{.2in}= \sum_{m,j} (-1)^{\maltese_j}\sF^{n-m+1}(a_1,\ldots,
  a_j, \mu_{\mA}^m(a_{j+1}. \ldots,a_{j+m}),a_{j+m+1},\ldots,a_n)
\end{align}
where $1 \leq m \leq n$ and $0 \leq j \leq n-m$. Given two $\ainf$
morphisms $\sF:\mA \to \mB$ and $\sG:\mA \to \mA$, the composition of
the two morphisms is given by the sequence of maps
\begin{equation} \label{ainfmorphcomp}(\sF^n \circ \sG^m)(\vec{a}) =
\sum_{j=0}^{n-1} (-1)^{(\maltese_j)||\sG^m||}\sF^n(1^{\tens j}\tens
\sG^m\tens 1^{\tens n-m-j})(\vec{a})\end{equation} for $\vec{a} \in
V^{\tens n+m-1}$ and $||\sG^m||$ the shifted degree of the map
$\sG^m$. We say that two $\ainf$ morphisms $\sF$ and $\sG$ between
$\mA$ and $\mB$ are \emph{$\ainf$ homotopic} if there is a sequence of
maps \[T^n: V^{\tens n} \to W[-n], \hspace{.2in} n \geq 1\] which
satisfy the sequence of relations
\begin{align}
\label{ainfhtpyrel}  &\sF^n(\vec{a})- \sG^n(\vec{a}) 
 = \sum_{1\leq r \leq n \atop{0\leq j \leq n-r}} (-1)^{\maltese_j}
 T^{n-r+1}(1^{\tens j}\tens \mu^r_{\mA}\tens 1^{\tens
   (n-r-j)})(\vec{a}) \\ \nonumber & \hspace{.2in} + \sum_{0\leq j,
   2\leq r \atop{j+r \leq n}}\sum_{S=n} (-1)^{\dagger}
 \mu^{j+r}_{\mB}(\sG^{s_1}\tens\ldots\tens\sG^{s_j}\tens
 T^{s_{j+1}}\tens \sF^{s_{j+2}}\tens
 \cdots\tens\sF^{s_{j+r}})(\vec{a})
\end{align}
where $\dagger = (|a_1| + \cdots + |a_{s_1+ \cdots + s_j}| - s_1 -
\cdots -s_j)$, $S=\sum_{k=1}^{j+r}s_k$, and $\vec{a} \in V^{\tens n}$.
These definitions follow the sign conventions of \cite{seidel08} for
$\ainf$ objects.

When there is no danger of confusion we shall omit the subscript on
the $\ainf$ composition maps. We shall always use $\mA$ and $\mB$ for
$\ainf$ algebras while $\sF$ and $\sG$ will always be $\ainf$
morphisms. Therefore $\sF^k$ will denote the $k$th term of $\sF$ and
hence be a map $\sF^k:V^{\tens k}\to V[1-k]$. Such a grading shift
means that given $a_1\tens \cdots \tens a_k \in V^{\tens k}$, where
$|a_i|$ denotes the grading of $a_i$, we have \[|\sF^k(a_1\tens \cdots
\tens a_k)| = \left(\sum_{i=1}^k|a_i|\right) + 1-k.\] For brevity, we
shall write multilinear combinations of multilinear maps with commas
rather than tensors and omit the input objects when working with
equalities.  We shall denote the differential graded algebra (dga) $V$
together with $\mu^1_{\mA}$ and $\mu^2_{\mA}$ as $A$. 

In Section~\ref{sec:families}, we define families of $\ainf$ algebras
indexed by a manifold $M$ as solutions to the Maurer-Cartan equation
on $\Omega^*(M;\g)$ where $\g$ is a particular locally trivial sheaf
of Lie algebras on $M$ built out of the Hochschild cochain complexes
of the underlying dgas of the $\ainf$ algebras. Theorem~\ref{thm:sFn}
gives an explicit formulation for the $\ainf$ morphism $\sF_{x\to
  y}:\mA_x \to \mA_y$ for $x,y \in M$ given a path $p:I\to M$
connecting $x$ and $y$. In Section~\ref{ahtpyonIxI}, after showing that differential homotopies
correspond to classical homotopies we prove Theorem~\ref{thm:hGn}
which gives an explicit form for a differential homotopy relating
families of $\ainf$ morphisms $\sF_t:\mA_x\to \mA_y$.

In Section~\ref{sec:examples}, we perform calculations. The goal is to
compute the cohomology of the total complex $\Omega^*(M;\g)$ in two
special cases. After a discussion of sheaf cohomology we argue that
for our purposes we can consider $BG$ as a manifold for $G$ finite or
finitely generated free nonabelian. Computation shows that in these
cases the spectral sequence collapses at $E^2$ as we discuss in
Section~\ref{ex:finitegroups} for a finite group and
Section~\ref{sec:nonabex} for a finitely generated free nonabelian
group. This leads up to Theorem~\ref{thm:alltrivial} which
states that for a finite group $\Gamma$, every homotopy
  $\Gamma$ action on an $\ainf$ algebra $\mA$ has class
  representatives $\sF_g:\mA\to \mA$ for all $g \in \Gamma$ which
  comprise a strict action. 

The author would like to thank Paul Seidel for his insightful questions and suggestions. 

\section{The differential graded Lie algebra $\Omega^*(M;\g)$}
\label{dglaOMg}
 \label{sec:families}
Let $V$ be a $\Z$-graded vector space and $M$ a differentiable manifold. Let
$\underline{V}$ denote a local system on $M$ with fibres isomorphic to
$V$.

Define \[\g_x = \prod_{n=0}^{\infty} \Hom (V_x^{\tens n},V_x[1-n]).\] Now
$\g_x$ is a $\Z$-graded Lie algebra under the Gerstenhaber bracket. Graded Lie
algebras are presented to good effect in \cite{gerstenhaber, goldmanmillson}. We have
chosen the indexing so that $\g_x^1$ contains the space
$\prod_{n=0}^{\infty}\Hom(V_x^{\tens n},V_x[2-n])$ of possible $\ainf$ structure
maps on $V_x$ including $\mu^0_x$. Let $\g$
denote the corresponding locally trivial sheaf of Lie algebras on $M$. Denote the algebra of $\g$ valued differential
forms on $M$ by $\Omega^*(M;\g)$.

\begin{proposition} \label{prop:Oisdgla} $\Omega^*(M;\g)$ is a differential $\Z$-graded Lie algebra with
  differential $d_{\nabla}$ and Lie bracket induced by the Gerstenhaber bracket
  on $\g$ and the wedge product of differential forms. \end{proposition}
\begin{proof}
First, because $\nabla$ is a flat connection we see that
$d_{\nabla}^2=0$ \cite[Appx C]{milnor}. Let $L^{m+k}= \Omega^m(M;\g^k)$ be the $(m+k)^{\text{th}}$ graded component for $m \geq
  0$ and $k \in \Z$. Then
  $d_{\nabla}(L^{m+k}) \subset \Omega^{m+1}(M;\g^k) \subset L^{m+k+1}$. Let
  $a=(\omega \tens \a) \in \Omega^m(M:\g^k)$ and $b=
  (\theta\tens \beta) \in \Omega^n(M;\g^{\ell})$, then   
  $d_{\nabla}[a,b] = [d_{\nabla}(a),b]+(-1)^{|a|}[a,d_{\nabla}(b)]$ by the
  Koszul rule of signs where $|a|$ denotes the total degree of $a$. Since the bracket is induced by the Gerstenhaber
  bracket on $\g$ together with the wedge product of differential forms we see
  that \begin{align*}\br{a}{b} &= \omega \wedge \theta \tens \br{\a}{\beta} \\
  &=  (-1)^{kn}\omega\wedge \theta \tens\a \circ \beta -
  (-1)^{|a||b|+m\ell}\theta\wedge\omega \tens \beta \circ \a \\
&  =  \omega \wedge \theta \tens \left((-1)^{nk}\a \circ \beta
  -(-1)^{k\ell+kn}\beta \circ \a\right)\in \Omega^{m+n}(M;\g^{k+\ell}).\end{align*} This shows that the bracket is linear
  with respect to total degree so $\br{L^i}{L^j} \subset L^{i+j}$. It is also
  clear that the bracket is homogeneous skew-symmetric because \begin{align*} -
  (-1)^{|a||b|}\br{a}{b} &= \omega\wedge \theta \tens \left( -(-1)^{k\ell+\ell m
  +mn}\a \circ \beta +(-1)^{\ell m+mn}\beta \circ \a\right) \\  &= \theta \wedge \omega \tens \left(-(-1)^{k\ell + \ell m}\a \circ \beta +
  (-1)^{\ell m}\beta \circ \a\right) \\
  &= \br{b}{a}n.\end{align*} Checking the Jacobi identity is a simple
  computation, so $\Omega^*(M;\g)$ is a dg Lie algebra
as desired.
\end{proof}

Consider a solution to the Maurer-Cartan equation
\begin{equation}\label{mc}
  d_{\nabla}(\a)+ \frac{1}{2}[\a,\a] = 0
  \end{equation} where $\a \in \Omega^m(M;\g^{1-m})$, or in other words the
  total degree of $\a$ is 1. Index the bigraded components
of $\a$ so that \[\a^{m,n} \in \Omega^m(M;\Hom(\underline{V}^{\tens n},\underline{V}[2-m-n])).\] Since $m
\geq 0$, this gives us a filter to use to understand solutions to
(\ref{mc}). When $m=0$, the $d_{\nabla}(\a)$ component cannot contribute, so
the $(0,n)$ component of (\ref{mc}) gives an element of $\Hom(V_x^{\tens
n},V_x[2-n])$ for each $x \in M$
so that:
\begin{align*} 0 &= \frac{1}{2}\sum_{m+r=n+1}[\a^{0,m},\a^{0,r}]\\ & =
\frac{1}{2}\sum_{m+r=n+1}\left(\sum_{j=0}^{m-1}(-1)^{\maltese_j}\a^{0,m}(1^{\tens
    j},\a^{0,r},1^{\tens n-r-j})\right. \\ &\hspace{1.4in}\left.+
  \sum_{k=0}^{r-1}(-1)^{\maltese_k}\a^{0,r}(1^{\tens k},\a^{0,m},1^{\tens
    n-r-m})\right) \\
&= \sum_{m+r=n+1}\sum_{j=0}^{m-1}(-1)^{\maltese_j}\a^{0,m}(1^{\tens
  j},\a^{0,r},1^{\tens n-j-r})
\end{align*}
which is precisely equation (\ref{ainfstruc}) on each point of $M$.
\begin{definition} 
A solution $\a \in \Omega^*(M;\g)$ to the Maurer-Cartan equation (\ref{mc})
with $|\a|=1$ gives a \emph{family of $\ainf$ algebras over $M$} where $\mA_x$
is the $\ainf$ algebra over $x$ for each $x \in M$ with the $\ainf$ structure
maps $\mu^n_x = \a^{0,n}_x$.\end{definition}

\begin{assumption} \label{an0is0} For ease of computation, we shall henceforth
assume that $\a^{0,0}=\a^{1,0}= \a^{2,0}=0$. This means that the curvature of each $\ainf$
algebra $\mA_x$ is zero as $\a^{0,0}=0$ identically. 
\end{assumption}

\begin{definition} For a manifold $M=K(G,1)$ with base point *, a family of $\ainf$ algebras over $M$
defines a \emph{homotopy group action by $G$} on the $\ainf$ algebra over the base
point where $[\sF_g]:\mA_* \to \mA_*$ for $g\in G$ is defined by integrating
around a loop corresponding to $g$ in $M$. \end{definition} We define such $\ainf$
morphisms in Theorem~\ref{thm:sFn}, and the homotopies between them in Theorem~\ref{thm:hGn}.  

\section{$\ainf$ morphisms on $I$}
\label{afunonI}

Let $\ga(t):[0,1] \to M$ be a path with $\ga(0)=x_0$ and $\ga(1)=x_1$. By
pulling back $\underline{V}$ along $\ga$ we may calculate on $I=[0,1]$. This is clear since
to determine $\sF:\mA_{x_0} \to \mA_{x_1}$ we shall integrate along $\ga$
and any tangent vectors perpendicular to $\ga$ will not contribute.
After pulling back to $I$, choose a trivialization compatible with the flat
connection $\nabla$. 

If we assume that $\a^{1,1}=0$, by integrating the first few levels of
the Maurer-Cartan equation (\ref{mc}) we calculate that 
\begin{align}
\nonumber \sF^1 & = \id, \\
\label{sfwde0} \sF^2 & = -\int_0^1\at^{1,2}, \\
\nonumber \sF^3 & = -\int_0^1\at^{1,3} - \int_{0\leq t \leq u \leq
  1}\a^{1,2}_u(\at^{1,2},1)+(-1)^{\maltese_1}\a^{1,2}_u(1,\at^{1,2}).
\end{align} 
 
From this point forward, we shall not include the signs in our
formulae as they all arise from the Gerstenhaber bracket and the
Koszul sign conventions.To prove a
general formula for the higher $\sF^n$ denote $\a^{1,n}$ as a height 1 tree
with $n$ leaves and a single root. We shall use the notation $(k,m)$-trees for
the sum of all height $k$ rooted trees with $m$ leaves where we do allow
valance 2 vertices corresponding to $\a^{1,1}$ terms.

\begin{definition}\label{sFnpd} Let $d$ be chosen so that
  $(\a^{1,1})^{d+1}=0$, then for $p\leq q \in [0,1]$:
  \begin{align*}
\sF^1_{p\to q} &=  (0,1)\text{-tree} + \sum_{i=1}^{d}\int_{p\leq
t_1 \leq \ldots \leq t_i \leq q}(i,1)\text{-trees} \\
\sF^n_{p\to q} &=\sum_{i=1}^{n-1+d} \int_{p\leq t_1 \leq \ldots \leq t_i \leq q}
(i,n)\text{-trees}. \end{align*} For $p > q$, simply reverse the inequalities in the integrals.\end{definition}

For example, when $d=0$ we have the initial terms given by (\ref{sfwde0}):
\[\sF^1 = \id, \hspace{1in} \sF^2 = \int_{1-\text{simplex}} \raisebox{-.15in}{\includegraphics[height=.4in]{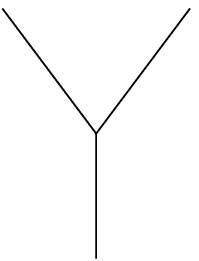}},\]
and
\[\sF^3 = \int_{1-\text{simplex}} \raisebox{-.15in}{\includegraphics[height=.4in]{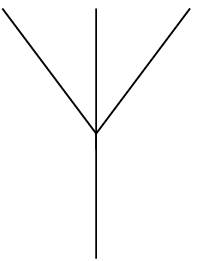}}
+ \int_{2-\text{simplex}}
\raisebox{-.15in}{\includegraphics[height=.4in]{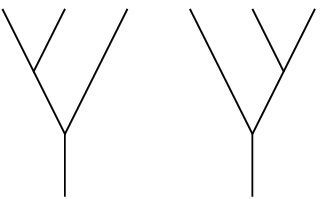}}.\]

\begin{lemma} \label{compositionpd}
For the $\sF^n$ defined in Definition~\ref{sFnpd}, composition follows the rule
\[\sF_{s\to t} \circ \sF_{r\to s} = \sF_{r\to t},\] for $s, r, t\in I$ where we
compose using (\ref{ainfmorphcomp}). 
\end{lemma}

\begin{proof} First, note that \[\sF^1_{r\to t} = \id + \sum_{i=1}^d\int_{r\leq
r_1 \leq \ldots \leq r_i \leq t}\hspace{-.5in}(i,1)\text{-trees},\] and
\begin{align*} \left(\sF_{s\to t}\right.&\left.\circ \sF_{r\to s}\right)^1 =  \sF_{s\to
t}^1\circ \sF_{r\to s}^1 \\
&= \left(\id + \sum_{j=1}^d\int_{s\leq
s_1 \leq \ldots \leq s_j \leq t}\hspace{-.5in}(j,1)\text{-trees}\right) \circ \left(\id + \sum_{i=1}^d\int_{r\leq
r_1 \leq \ldots \leq r_i \leq s}\hspace{-.5in}(i,1)\text{-trees}\right)  \\
& = \id_{r\to t} + \id_{s\to t} \circ \sum_{i=1}^d\int_{r\leq
r_1 \leq \ldots \leq r_i \leq s}\hspace{-.6in}(i,1)\text{-trees} + \left(\sum_{i=1}^d\int_{s\leq
s_1 \leq \ldots \leq s_i \leq t}\hspace{-.5in}(i,1)\text{-trees}\right)\circ \id_{r\to s} \\
&\hspace{.15in} + \sum_{\substack{j+i=2 \\ j \geq 1, i\geq 1}}^{2d} \int_{r \leq
r_1 \leq \ldots \leq r_i \leq s \leq s_1 \leq \ldots \leq s_j \leq
t}\hspace{-1.2in}(i+j,1)\text{-trees} \\
&= \sF^1_{r\to t}
\end{align*}
because the last sum is zero for $j+i > d$. We shall use a similar approach for
  higher $n$ and split up the desired result into pieces with mixed components
  or terms in only one half or the other.  The $n$th term of $\sF_{r\to t}$
  consists of integrals over the appropriate simplices of all $n$-leaved
  trees. On the other side we have
  \begin{align*}(&\sF_{s\to t} \circ \sF_{r\to
  s})^n = \sum_{i=1}^n \sum_{(\sum m_j)=n}\sF_{s\to t}^i(\sF_{r\to
  s}^{m_1}\tens \cdots \tens \sF_{r\to s}^{m_i})\\
    &= \sum_{i=1}^n \sum_{(\sum m_{\ell})=n} \left(\sum_{k=1}^{i-1+d} \int_{s\le
  s_1 \le \cdots \le s_k \le
  t} \hspace{-.5in}(k,i)\text{-trees}\right) \circ \\
  &\hspace{1in}\left(\sum_{j_1=1}^{m_1-1+d} \int_{r\le r_1 \le \cdots \le
  r_{j_1} \le s}\hspace{-.5in}
  (j_1,m_1)\text{-trees} \tens \cdots \tens \sum_{j_i=1}^{m_i-1+d} \int_{r\le
  r_1 \le \cdots \le r_{j_i} \le s}\hspace{-.5in}
  (j_i,m_i)\text{-trees}\right) \\
  &\hspace{.3in}+ \id_{s\to t}(\sF_{r\to
  s}^{n}) + \sF_{s\to t}^n(\id_{r\to
  s}\tens \cdots \tens \id_{r\to s}) .  \end{align*}
As in the $n=1$ example, we
  shall need to include identities to transport values to the appropriate endpoints. This allows us to
  increase the depth of a tree freely. Now, let $j$ be the maximum depth of the
  $\{j_i\}$ of any fixed configuration. Extend the other trees with identities
  to the same depth so we can combine the terms. This is a finite process
  because we do not modify the maximum-depth tree so the depth of the combined
  tree is simply $j+k$. We shall also use the fact that any tree deeper than
  $n-1+d$ is automatically 0. Thus, we can continue the above equation as
  follows:
  \begin{align*}  &= \sum_{i=1}^n \sum_{(\sum
  m_{\ell})=n} \sum_{\substack{j+k =2,\\ j,k\ge 1}}^{n-1+d} \int_{s\le
  s_1 \le \cdots \le s_k \le
  t} \hspace{-.5in}(k,i)\text{-trees}\left(\int_{r\le r_1 \le \cdots \le
  r_{j} \le s}\hspace{-.7in} (j,m_1)\text{-trees} \tens \cdots \tens \int_{r\le
  r_1 \le \cdots \le r_{j} \le s}\hspace{-.7in} (j,m_i)\text{-trees}\right)\\
  &\hspace{.2in} + \sum_{k=1}^{n-1+d}\int_{r\le r_1 \le \cdots \le r_k \le
  s} \hspace{-.5in}(k,n)\text{-trees} + \sum_{k=1}^{n-1+d}\int_{s\le
  s_1 \le \cdots \le s_k \le t} \hspace{-.5in}(k,n)\text{-trees} \\
  &= \sum_{k=1}^{n-1+d}\sum_{\ell = 0}^k \int_{r=r_0 \le r_1 \le \cdots \le
  r_{\ell} \le s \le r_{\ell +1} \cdots \le r_k \le t} \hspace{-1.5in}
  (k,n)\text{-trees}\\ &= \sum_{k=1}^{n-1+d} \int_{r\le r_1 \le \cdots \le
  r_k \le t} \hspace{-.5in}(k,n)\text{-trees} \\ &= \sF^n_{r\to t} \end{align*}
\end{proof}
\begin{corollary}\label{cor:strictinverses} Given a path $p:I\to M$, the strict inverse of $\sF_{p(t)}$ is
  $\sF_{p(1-t)}$ because $\sF_{r\to s}\circ \sF_{s\to r}=\sF_{s\to
  s}=\id$. \end{corollary}

\begin{theorem} \label{thm:sFn} When $\a^{1,1}=0$, the maps $\sF^k$ defined in
Definition~\ref{sFnpd} are $\ainf$ morphisms from $\mA_p$ to $\mA_q$.
\end{theorem}

Lemma~\ref{compositionpd} demonstrates that the proposed morphism maps in
Theorem~\ref{thm:sFn} compose as required for $\ainf$ morphisms in (\ref{ainfmorphcomp}). Using this
composition property, we shall apply results from differential equations to argue
that our choice is correct.


\begin{lemma}\label{diffeq}Using the $\sF^k$ in Theorem~\ref{thm:sFn}, define \begin{align*}\Delta_{r,s}^n := \sum_{i=2}^{n}\sum_{\sum k_j =n}& \mu^i_s\left(\sF^{k_1}_{r,s} \tens \cdots \tens
\sF^{k_i}_{r,s}\right) \\ &\hspace{.8in}- \sum_{i=2}^{n}\sum_{j=0}^{n-i} (-1)^{\maltese_j}\sF^{n-i+1}_{r,s}\left(1^{\tens j}\tens \mu^{i}_r
\tens 1^{\tens n-(i+j)}\right).\end{align*}
Then \[\Delta_{r,s}^n = 0 ~\forall s,n.\]
\end{lemma}

Notice that $\Delta_{r,s}^n$ is precisely a measurement of the failure of our
proposed maps to satisfy the $n^{\text{th}}$ level of
the $\ainf$ morphism relations (\ref{ainfmorphrel}). 
\begin{proof}
  First, we determine the initial conditions. Clearly, $\Delta_{r,r}^n = 0$ as
  $\sF^1=\id$ and all higher order terms $\sF^k$ are zero because they involve integrating over sets
  of measure zero. Now
\begin{align} \label{dDelta0s}
\nonumber \left.\der{}{s}\Delta_{r,s}^n\right|_{s=r} &=
 \left.\der{\mu_s^n}{s}\right|_{s=r} +
 \sum_{k=2}^{n-1}\sum_{j=0}^{n-k}\left.(-1)^{\maltese_j}\mu_r^{n-k+1}\left(1^{\tens
   j}\tens \der{\sF^k_{r\to s}}{s}\right|_{s=r} \tens 1^{\tens n-j-k}\right)
 \\ &\hspace{.15in}-
 \sum_{k=2}^{n-1}\sum_{j=0}^{k-1}\left.(-1)^{\maltese_j}\der{\sF^k_{r\to s}}{s}\right|_{s=r}\left(1^{\tens
   j},\mu_r^{n-k+1}, 1^{\tens k-j-1}\right) \\
\nonumber &= \left.\der{\mu_s^n}{s}\right|_{s=r} + \sum_{k=2}^{n-1}\br{\a^{0,n-k+1}_r}{\a_r^{1,k}} = 0
\end{align}
where we used (\ref{mc}) at the point $x=r$ on $I$ and the fact that
\begin{equation} \label{dsF0s}\left.\der{\sF^k_{r\to
    s}}{s}\right|_{s=r} = \begin{cases}\left.\der{}{s}\left(\int_r^s\a^{1,*}_t dt +
\mathcal{O}((s-r)^2)\right)\right|_{s=r} =\a_r^{1,*} & \text{if } k\geq 2 \\
\left.\der{}{s}\id\right|_{s=r} = 0 & \text{if } k=1. \end{cases}
\end{equation}
Determine a differential equation defining $\Delta_{r,s}^k$ in terms of
$\Delta_{0,s}$. 
\begin{align*}
\Delta_{0,s}^n &= \sum_{i=2}^{n}\sum_{K_i =n}
\mu^i_s\left(\sF^{k_1}_{0,s} \tens \cdots \tens \sF^{k_i}_{0,s}\right) \\
& \hspace{.55in}- \sum_{i=2}^{n}\sum_{j=0}^{n-i} (-1)^{\maltese_j}\sF^{n-i+1}_{0,s}\left(1^{\tens j}\tens \mu^{i}_0
\tens 1^{\tens n-i-j}\right) \\
&= \sum_{i=2}^n\sum_{K_i =n}
\mu^i_s\left(\sum_{M_{\ell} = k_1}\hspace{-.07in}\sF^{\ell}_{r,s}(\sF_{0,
  r}^{m_1}\tens\cdots\tens\sF_{0,r}^{m_{\ell}}) \tens
\cdots\right. \\ &\hspace{2.5in}\left.\cdots \tens \hspace{-.07in}\sum_{M_{\ell}=k_i}\hspace{-.07in}\sF^{\ell}_{r,s}(\sF_{0,
  r}^{m_1}\tens\cdots\tens\sF_{0,r}^{m_{\ell}})\right) \\
&\hspace{.15in} -
\sum_{i=2}^n\sum_{M_k=n-i+1}\sum_{\ell=1}^{k}\sum_{j=0}^{m_{\ell}-1}(-1)^{\maltese_{M_{\ell-1}+j}}
\sF^{k}_{r,s} \circ \\ & \hspace{.8in}\left(\sF^{m_1}_{0,r}\tens\cdots\tens\sF^{m_{\ell}}_{0,r}(1^{\tens
  j}, \mu^{i}_0,1^{\tens m_{\ell}-1-j})\tens \cdots 
\tens \sF^{m_k}_{0,r}\right)
\end{align*}
where we use the notation $M_k = \sum_{i=1}^k m_i$. Therefore,
\begin{align*}
&\Delta^n_{0,s}
= \sum_{k=1}^n\Delta^k_{r,s}\left(\sum_{M_k = n}(\sF_{0,
  r}^{m_1}\tens\cdots\tens\sF_{0,r}^{m_k})\right) \\ &\hspace{.05in}+ 
\sum_{k=1}^n\sum_{M_k=n} \sum_{\ell=1}^k \sF^{k}_{r,s}\left(\sF^{m_1}_{0,r}\tens\cdots\tens\sF^{m_{\ell-1}}_{0,r}\tens \Delta_{0,r}^{m_{\ell}}\tens \sF^{m_{\ell+1}}_{0,r}\tens \cdots 
\tens \sF^{m_j}_{0,r}\right).
\end{align*}
Thus
 \begin{align*}
 &\der{\Delta^n_{0,s}}{s}
 =  \sum_{k=1}^n\der{\Delta^k_{r,s}}{s}\left(\sum_{M_k=n}
 (\sF_{0,r}^{m_1}\tens\cdots\tens\sF_{0,r}^{m_k})\right)\\
 &\hspace{.1in}+ \sum_{k=1}^n\sum_{M_k=n} \sum_{\ell=1}^k \der{\sF^{k}_{r,s}}{s}\left(\sF^{m_1}_{0,r}\tens\cdots\tens\sF^{m_{\ell-1}}_{0,r}\tens \Delta_{0,r}^{m_{\ell}}\tens \sF^{m_{\ell+1}}_{0,r}\tens \cdots 
\tens \sF^{m_j}_{0,r}\right)
 \end{align*}
 since all the other terms are independent of $s$ and hence simply
 disappear. Now $r$ is simply a variable so we can choose $r=s$. After applying
 (\ref{dDelta0s}) and (\ref{dsF0s}) we have 
 \begin{align*}
 \der{\Delta^n_{0,s}}{s}
 &= \sum_{k=2}^n\sum_{M_k=n} \sum_{\ell=1}^k \a^{1,k}_s\left(\sF^{m_1}_{0,s}\tens\cdots\tens\sF^{m_{\ell-1}}_{0,s}\tens \Delta_{0,s}^{m_{\ell}}\tens \sF^{m_{\ell+1}}_{0,s}\tens \cdots 
\tens \sF^{m_j}_{0,s}\right).
 \end{align*}
Since $k \geq 2$, we see that $m_{\ell}<n$. We also know by explicit
computation that $\Delta^2_{r,s}=0$ and $\Delta^3_{r,s}=0$. By the
multilinearity of $\a^{1,k}_s$, we therefore know that
$\der{\Delta^n_{0,s}}{s}=0$ for $n=2,3,4$. Now we can induct on $n$. Say we
know that $\Delta^k_{0,s}=0$ for all $k <n$ and hence that
$\der{\Delta^k_{0,s}}{s}=0$ for $k \leq n$. Using the initial condition that
$\Delta^n_{0,0}=0$, by the existence and uniqueness of linear ordinary
differential equations, it is clear
that $\Delta^n_{0,s}=0$ for all $s$. However, there was nothing special about
our choice of 0 as the initial point so by translating by $r$, 
 $\Delta_{r,s} = 0$ for all $r,s \in I$. \end{proof}

 \begin{proof}[Proof of Theorem~\ref{thm:sFn}] By Lemma \ref{compositionpd},
  the manipulations that we took to prove Lemma \ref{diffeq} were
 valid. Therefore, $\sF_{r,s}$ is an
 $\ainf$ morphism from $\mA^*_p$ to $\mA^*_q$.
\end{proof}

Note that while we assumed $\a^{1,1}=0$ in Theorem~\ref{thm:sFn}, we
did so to eliminate the possibility of infinite recursions with
Stokes' Theorem leading to infinite depth trees with valence two
vertices in Definition~\ref{sFnpd}. We may relax the assumption that
$\a^{1,1}=0$ by instead requiring a finite descending descending
filtration on $\Omega^*(M;\g)$ with the condition that $\a^{1,1}$, and
any non-invariant component of $\a^{0,1}$, decrease the degree. This
will give us a nilpotency restriction on $\a^{1,1}$ thereby fixing a
$d<\infty$ for Definition~\ref{sFnpd}.

\section{$\ainf$ homotopies on $M=I\times I$} \label{ahtpyonIxI}

\subsection{Classical homotopies}
Following \cite[Ch. X]{grifmorg}, \cite[Ch. I]{seidel08} we can algebraically define a (classical) homotopy between $\sF_0:\mA\to \mB$
and $\sF_1:\mA\to \mB$ to be given by the strictly commutative diagram
\[\xymatrix{
&  \mB \\ \mA \ar@/^1pc/[ru]^{\sF_0} \ar@/_1pc/[rd]_{\sF_1} \ar[r]& I \tens \mB \ar[u]_{ev_0}
  \ar[d]^{ev_1}  \\  
&\mB}\]
where $I$ is considered as the quiver
\[\xymatrix{
^{u_0}\bullet \ar[r]^{h} & \bullet^{u_1}
}\]
with the relations \begin{align}
&|u_0| = |u_1| = 0, \hspace{.2in} |h|=1 \nonumber \\
& u_0^2 = u_0 \hspace{.2in} u_1^2 = u_1 \hspace{.2in} h^2 = 0 \label{Imult}\\
& hu_0 = h \hspace{.2in} u_0h = 0 = hu_1 \hspace{.2in} u_1h = h. \nonumber 
\end{align}
We then impose a differential $\partial$ according to the rules:
\begin{align} \label{Idiff} \partial u_0 &= h = -\partial u_1 \hspace{.2in} \partial h = 0.
\end{align}

Note that $I$ is a differential graded algebra (dga). This is clear by
  checking the interaction of $\partial$ with the relations in equations
  \ref{Imult} and \ref{Idiff} as well as the Leibnitz rule.

\subsection{Differential homotopies} \label{sec:diffhtpy}
By extrapolating the relevant characteristics of the classical picture, we can
define differential homotopies on a family of $\ainf$ morphisms between $\mA$
and $\mB$ indexed by $t \in [0,1]$ by the commutative diagram of $\ainf$
algebras below:
\[\xymatrix{ & \mB \\ \mA \ar@/^1pc/[ru]
\ar@/_1pc/[dr] \ar[r]^-{\phi} & \Omega^*([0,1],\mB). \ar[u] \ar[d] \\ &\mB}\]
We shall sometimes denote $\Omega^*([0,1],\mB)$ by
$\Omega^*([0,1])\tens \mB$ in the sense of \cite{marklshnider}
or \cite{seidel08} since $\Omega^*([0,1])$ is a dga. The $\ainf$ structure maps
on the tensor product are given by the formulae:
\begin{align*} \mu^1_{O\mB}(a\tens b) &:= d(a)\tens b +
(-1)^{|a|} a\tens\mu^1_{\mB}(b)  \\
\mu^n_{O\mB}(a_1\tens b_1, \cdots, a_n\tens b_n) &:= (-1)^{\diamond}a_1\cdots a_n \tens
\mu_{\mB}^n(b_1,\cdots, b_n) \hspace{.2in} n\geq 2 
\end{align*}
where $\diamond = \sum_{j<k} |b_j||a_k|$
since $\Omega^*([0,1])$ is a dga. Notice that $\diamond \neq 0$ only if exactly
one $a_k \in \Omega^1([0,1])$. \begin{definition}Let $\sF_t:\mA \to \mB$ be our family of
$\ainf$ morphisms with \[\phi^n(a_1,\cdots,a_n)(t) = \sF_t^n(a_1,\cdots,a_n) +
    \Theta_t^n(a_1,\cdots,a_n)dt.\] Then the $\Theta_t^n$, $n \geq 1$ form a \emph{differential homotopy} with respect to the family $\sF_t$. \end{definition} \begin{remark} Now, since $\phi$ is an $\ainf$
homomorphism, by (\ref{ainfmorphrel}) we have, for $n \geq 1$ the relations
\begin{align*}
\sum_{k=1}^n\sum_{\sum r_i=n}\mu_{O\mB}^k(\phi^{r_1},\cdots,\phi^{r_k}) &=
\sum_{i=1}^n\sum_{j=0}^{n-i} \phi^{n-i+1}(1^{\tens j},\mu_{\mA}^i,1^{\tens
  n-i-j} )
\end{align*}
which become
\begin{align*}
&\sum_{k=1}^n\sum_{\sum r_i=n}\mu^k_{O\mB}(\sF_t^{r_1} + \Theta_t^{r_1}dt,\cdots,\sF_t^{r_i} + \Theta_t^{r_k}dt) \\&\hspace{.3in}=
\sum_{i=1}^n\sum_{j=0}^{n-i}\sF_t^{n-i+1}(1^{\tens
  j},\mu_{\mA}^i, 1^{\tens n-i-j}) 
 + \sum_{i=1}^n\sum_{j=0}^{n-i}\Theta_t^{n-i+1}(1^{\tens
  j},\mu_{\mA}^i, 1^{\tens n-i-j})dt. 
\end{align*}
Applying the definition of $\mu_{O\mB}^k$ we see that 
in degree 0 each $\sF_t$ is required to be an $\ainf$ morphism
between $\mA$ and $\mB$. However, in the coefficients of $dt$ we have
\begin{align*}\label{eqn:diffhtpy}
\der{\sF_t^n}{t} +
\mu_{\mB}^1(\Theta_t^n)&= \sum_{i=1}^n\sum_{j=0}^{n-i}\Theta_t^{n-i+1}(1^{\tens
  j},\mu_{\mA}^i,1^{\tens n-i-j}) \\ \nonumber &\hspace{.15in} -\sum_{k=1}^n\sum_{\sum r_i=n}\sum_{\ell=1}^k
\mu_{\mB}^k(\sF_t^{r_1},\cdots,\sF_t^{r_{\ell-1}},\Theta_t^{r_{\ell}},\sF_t^{r_{\ell}+1},\cdots,\sF_t^{r_k}).
\end{align*}
This is the differential equation is quite similar to the standard $\ainf$ homotopy relation given in (\ref{ainfhtpyrel}). \end{remark} 

\begin{proposition}\label{prop:diffhtpytoclassical} Differential homotopies give rise to classical homotopies.
  \end{proposition}
\begin{proof}
We must show strict commutativity of the following diagram
\[\xymatrix{
\R& I\ar[l]_{\text{ev}_0} \ar[r]^{\text{ev}_1} & \R\\
& \Omega^*([0,1])\ar@/^1pc/[lu]\ar@/_1pc/[ru]\ar@.[u]_{q.i.}^{\Phi} &} \]
for a quasi-isomorphism $\Phi$. We shall denote
$\Omega^*([0,1])$ by $\mA$ for ease of notation. 
We proceed by defining $\Phi^n$ inductively as we require $\Phi$ to be
an $\ainf$ morphism satisfying (\ref{ainfmorphrel}). For $a, b \in \mA^0=\Omega^0([0,1])$, define $\Phi^1$
by \[\Phi^1(a+bdt) = a(0)u_0+a(1)u_1 + \left(\int_0^1bdt\right)h.\]
Thus it is clear that $\Phi^1(fg)=\Phi^1(f)\Phi^1(g)$ for $f,g \in \mA^j$. Therefore,
we need to have a $\Phi^2$ that will cancel the mixed terms
\[\Phi^1(a(t)b(t)dt) = \left(\int_0^1 a(t)b(t)dt\right) h \] and \[\Phi^1(a(t))\Phi^1(b(t)dt) = \left(a(1)\int_0^1b(t)dt\right)h\]
exactly without disturbing the equality of our earlier compositions. 
Recall that $\mu^2_{\mA}$ is just normal multiplication of forms and
$\mu^1_{\mA}$ is differentiation of forms. Also, $\mu^2_{I}$ is the
linear extension of our multiplication table given in equation
\ref{Imult} and $\mu^1_{I} = \partial$ as constructed in equation
\ref{Idiff}. Hence, the only possible input for which $\Phi^2$ should
be nonzero is $(f(t)dt, g(t)dt)$ as all others would contribute to
equations that are already satisfied. Let \begin{equation}\label{eqn:phi2}\Phi^2(f(t)dt,g(t)dt) :=
\left(\int_{0\leq t\leq s \leq 1} f(s)g(t)ds dt\right)h.\end{equation} Now,
putting all of this together, we require
\begin{align*}
\Phi^1(\mu^2_{\mA})(f(t),g(t)dt) &+ \Phi^2(\mu^1_{\mA}\tens
1+1\tens\mu^1_{\mA})(f(t),g(t)dt) \\ &=
\mu_{I}^1(\Phi^2)(f(t),g(t)dt) +
\mu^2_{I}(\Phi^1\tens\Phi^1)(f(t),g(t)dt)
\end{align*}
or in other words,
\begin{align*}
\left(\int_0^1f(t)g(t)dt\right)h &+ \Phi^2\left(\der{f(t)}{t}dt,g(t)dt\right) +0 = 0 +
\left(f(1)\int_0^1 g(t)dt\right)h. 
\end{align*}
However, by Stokes' Theorem and (\ref{eqn:phi2}), this is satisfied. The check
that \[\Phi^1(g(t)dtf(t))+\Phi^1(g(t)dt)\Phi^1(f(t)) = \Phi^2(1\tens
  \mu^1_{\mA}+\mu^1_{\mA}\tens 1)(g(t)dt,f(t))\] is similar.
 Thus, we may define
 \begin{align} \label{genphi} \Phi^k(&a_1(t),\cdots,a_k(t))  =\nonumber \\ &=\begin{cases}
\left(\int_{0\leq t_k \leq \cdots \leq t_1\leq 1} a_1(t_1)\cdots
a_k(t_k)\right)h & \text{if~}a_i(t) \in \mA^1 ~\forall i \\
0 & \text{if any~} a_i \in \mA^0
\end{cases}.\end{align}
 
In the $\ainf$ morphism relations (\ref{ainfmorphrel}), all terms involving $\mu^n, n\geq 3$
drop out, and the term $\mu^1_{I}(\Phi^{k+1})$ drops out by $\partial h =
0$ any time $\Phi^{k+1}$ is nonzero. Hence, 
the only terms which are nonzero on the left (LHS)
have exactly one 0-form entry and the only nonzero terms on
the right (RHS) are those involving $\Phi^1$ and $\Phi^k$. If we have more
than one 0-form, then the LHS is also zero because of our definitions
of any $\Phi^n, n>1$ since we have explicitly checked $\Phi^2$. On the
other hand, if there are no 0-form entries, then the LHS is zero and
the RHS is $O(h^2)=0$. Therefore, the relations reduce to:
\begin{align*}\Phi^{k+1}\left(\sum_{j=0}^k 1^{\tens j}\tens\mu^1_{\mA}\tens 1^{k-j}\right) &+
\Phi^k\left(\sum_{j=0}^{k-1} 1^{\tens j}\tens\mu^2_{\mA}\tens 1^{k-1-j}\right) \\&= \mu^2_{I}(\Phi^1\tens \Phi^k + \Phi^k\tens\Phi^1)\end{align*}
where, for $f_j \in \mA^0$, the possible arguments are
\[(f_1(t)dt,\cdots,f_{i-1}(t)dt,f_i(t),f_{i+1}(t)dt,\cdots,f_{k+1}(t)dt),
\hspace{.2in} i=1,\cdots,k+1.\]

Thus, on the left we have
\begin{align*}
&LHS = \left(\int_{0\leq t_{k+1}\leq\cdots\leq t_{i+1}=t_i\leq \cdots \leq t_1\leq
  1}\hspace{-.8in}f_1(t_1)\cdots f_i(t_i)\cdots f_{k+1}(t_{k+1})dt_1\cdots
  dt_{i-1}dt_{i+1}\cdots dt_{k+1}\right.\\
&\hspace{.15in} +\int_{0\leq t_{k+1}\leq\cdots\leq t_{i}=t_{i-1}\leq \cdots \leq t_1\leq
  1}\hspace{-.8in}f_1(t_1)\cdots f_i(t_i)\cdots f_{k+1}(t_{k+1})dt_1\cdots
  dt_{i-1}dt_{i+1}\cdots dt_{k+1}\\
&\hspace{.15in} + \int_{0\leq t_{k}\leq\cdots\leq t_1\leq
  1}f_1(t_1)\cdots f_{i-1}(t_{i-1})f_i(t_{i-1})f_{i+1}(t_i)\cdots f_{k+1}(t_{k})dt_1\cdots
dt_{k}\\
&\left.\hspace{.15in} + \int_{0\leq t_{k}\leq\cdots\leq t_1\leq
  1}f_1(t_1)\cdots f_i(t_{i})f_{i+1}(t_i)f_{i+2}(t_{i+1})\cdots f_{k+1}(t_{k})dt_1\cdots
dt_{k}\right)h
\end{align*}
which agrees with the RHS in all cases.
\end{proof}

\begin{remark} Since $\Phi:\Omega^*([0,1]) \to I$ is a quasi-isomorphism, it is
true that differential homotopies correspond to classical homotopies by
inverting $\Phi$, but we shall not do the necessary calculations here.
\end{remark}

\subsection{Mapping from $I\times I $ to $I$}
\begin{theorem}\label{IItoI}
  Let $\a^{*,*}$ be a solution of equation \ref{mc} on $I\times I$. We can
  collapse $I\times I$ to $I$ by applying the standard
  projection onto the first element. We shall identify the first interval with
  $0\leq s \leq 1$ and the second with $0 \leq t \leq 1$. Then we
  can construct a solution $\hat{\a}^{m,n} \in
  \Omega^m(I,\Hom((\Omega^*(I,V))^{\tens n},\Omega^*(I,V))$ on the
  submanifold $I$ as follows:
\begin{align*}
\hat{\a}^{0,n} &= \delta_{1,n}d + \a^{0,n} + \iota_{\partial_t}\a^{1,n}dt \\
\hat{\a}^{1,n} &= \iota_{\partial_s}\a^{1,n}ds + (1\tens
\iota_{\partial_t}+\iota_{\partial_t}\tens 1)\a^{2,n}dt \\
\hat{\a}^{m,n} &= 0 \hspace{.2in} \text{for all}~m \geq 2
\end{align*}
where $d$ is the differential on $\Omega^*(I,V)$ considered as a
constant function of $s$, and $\delta$ is the
Kronecker delta function.
\end{theorem}
\begin{proof}
 We must check that $\hat{\a}^{k,\ell}$ is a solution to (\ref{mc}). Let us begin
 with the lowest level when $k=0$.
\begin{align*}
 \sum_{m+r=n+1}&\frac{1}{2}\left[\hat{\a}^{0,m},\hat{\a}^{0,r}\right]
\\
 &= \sum_{m+r=n+1}\frac{1}{2}\left(\left[ \delta_{1,m}d,\delta_{1,r}d\right] + \left[\a^{0,m},\delta_{1,r}d\right]+
  \left[\iota_{\partial_t}\a^{1,m}dt, \delta_{1,r}d\right] 
  \right.\\ &\hspace{.6in}+ [\delta_{1,m}d,\iota_{\partial_t}\a^{1,r}dt]+ \left[\delta_{1,m}d,\a^{0,r}\right]
  +\left[\a^{0,m},\a^{0,r}\right]  \\ & \hspace{.6in}+
  \left[\a^{0,m},\iota_{\partial_t}\a^{1,r}dt\right]+\left.\left[\iota_{\partial_t}\a^{1,m}dt,\a^{0,r}\right]
  +\left[\iota_{\partial_t}\a^{1,m}dt,\iota_{\partial_t}\a^{1,r}dt\right]\right)\\
&= \sum_{m+r=n+1}\left( \left[\delta_{1,m}d,\a^{0,r}\right]
   +\left[\iota_{\partial_t}\a^{1,m}dt,\a^{0,r}\right]\right)\\
&= \iota_{\partial_t}\left(d_{\nabla}\a^{0,n}+ \sum_{m+r=n+1}\br{\a^{1,m}}{\a^{0,r}}\right)dt \\
&=0.
\end{align*}
Now consider the component of (\ref{mc}) that lies in $\Omega^1(I)$.
\begin{align*}
\der{}{s}(\hat{\a}^{0,n})ds &+ \sum_{r+m=n+1}[\hat{\a}^{0,m},\hat{\a}^{1,r}] \\
&=\der{(\a^{0,n}+\iota_{\partial_t}\a^{1,n})}{s}ds + \der{(\iota_{\partial_s}\a^{1,n})}{t}dsdt + \sum_{r+m=n+1}\left([\a^{0,m},\iota_{\partial_s}\a^{1,r}ds]\right. \\ & \left. +    [\iota_{\partial_t}\a^{1,m}dt,\iota_{\partial_s}\a^{1,r}ds] + + [\a^{0,m},(1\tens
  \iota_{\partial_t}+\iota_{\partial_t}\tens 1)\a^{2,r}dt]\right) \\
& = \iota_{\partial_s}\left(d_{\nabla} \a^{0,n} + \sum_{m+r=n+1} [\a^{0,m},\a^{1,r}]\right)ds  \\
& + (1\tens
  \iota_{\partial_t}+\iota_{\partial_t}\tens 1)\left(d_{\nabla}\a^{1,n} + \sum_{m+r=n+1}[\a^{0,m},\a^{2,r}]+\frac{1}{2}[\a^{1,m},\a^{1,r}]\right) \\
&=0.
\end{align*}
All further levels of the Maurer-Cartan equation are zero by dimensionality.
\end{proof}

\subsection{$\ainf$ homotopies}
Let $\a^{m,n}$ be a solution to the Maurer-Cartan equation on the
square $I\times I$ with $\a^{1,1}=0$. Since we presume this square is
a pullback of two homotopic paths on $M$, we require that the
$\ainf$ structure be constant for the edges $s=0$ and $s=1$.  Let
$\sF_t$ be a family of $\ainf$ morphisms given by
Theorem~\ref{thm:sFn} using paths along constant $t$ between the
$\ainf$ algebras $\mA_{s=0}$ and $\mA_{s=1}$ on the square. Now $\hat{\a}^{0,1}= d$ and $\hat{\a}^{1,1}=(1\tens
\iota_{\partial_t}+\iota_{\partial_t}\tens 1)\a^{2,1}$ are both nilpotent of order 2.

\begin{theorem} \label{thm:hGn} Define $\hG^n$ using the $\ha^{1,k}$ terms in
Definition~\ref{sFnpd} with $d=1$. Then
  $\hat{\sG}$ defines a differential homotopy with respect to the $\sF_t$.
  \end{theorem}
\begin{proof}
First, we show that $\hG$ is an $\ainf$ morphism from
$\Omega^*(I;\mA_{s=0})$ to $\Omega^*(I;\mA_{s=1})$. Then we shall show that
there is a natural injection of $\mA_{s=0}$ into $\Omega^*(I;\mA_{s=0})$.

By definition, $\hG^n = \sF^n + \sG^n$ where
$\text{Im}(\sG^n) \in \Omega^1(I;\mA_{s=1})$. As in Lemma~\ref{diffeq},
define \begin{align}\nonumber \hat{\Delta}_{r,s}^n &:= \sum_{i=2}^{n}\sum_{\sum k_j =n} \mu^i_s\left(\hG^{k_1}_{r,s} \tens \cdots \tens
\hG^{k_i}_{r,s}\right) \\ \nonumber &\hspace{1in}- \sum_{i=2}^{n}\sum_{j=0}^{n-i} (-1)^{\maltese_j}\hG^{n-i+1}_{r,s}\left(1^{\tens j}\tens \mu^{i}_r
\tens 1^{\tens n-(i+j)}\right)\\
\label{hDeltaindt}&= \Delta^n_{r,s}+\sum_{i=2}^{n}\sum_{\ell=1}^i\sum_{\sum k_j
=n} \mu^i_s\left(\sF^{k_1}_{r,s} \tens \cdots \tens \sG^{k_{\ell}}_{r,s} \tens
\cdots \tens \sF^{k_i}_{r,s}\right) \\ \nonumber  
&\hspace{1in}- \sum_{i=2}^{n}\sum_{j=0}^{n-i}
(-1)^{\maltese_j}\sG^{n-i+1}_{r,s}\left(1^{\tens j}\tens \mu^{i}_r 
\tens 1^{\tens n-(i+j)}\right).
\end{align}
Since $\Delta^n_{r,s}=0$ by Lemma~\ref{diffeq}, we see that
$\hat{\Delta}^n_{r,s} \in \Omega^1(I;\mA_{s=1})$. However in this case,
\[\left.\der{\hG^k_{r,s}}{s}\right|_{s=r} = \ha^{1,k}_r,\] so we also have
\begin{align*}\der{\hat{\Delta}^n_{0,s}}{s}
&= \sum_{k=1}^n\sum_{M_k=n}\sum_{\ell=1}^k \ha^{1,k}_s(\hG^{m_1}_{0,s}\tens \cdots \tens \hG^{m_{\ell-1}}_{0,s}\tens \hat{\Delta}^{m_{\ell}}_{0,s}\tens \hG^{m_{\ell+1}}_{0,s}\tens \cdots \tens \hG^{m_j}_{0,s})
\end{align*}
As before, proceed by induction. We know that
\[\hat{\Delta}^1_{0,s}= \ha^{0,1}_s(\sG^1_{r,s})-\sG^1_{r,s}(\ha^{0,1}_r)=0.\]
Assume that $\hat{\Delta}^k_{0,s}=0$ for $k < n$. For $k=n$, after using
(\ref{hDeltaindt}) we are thus left with
\[\der{\hat{\Delta}^n_{0,s}}{s}
=  \ha^{1,1}_s(\hat{\Delta}^{n}_{0,s})
=0.\] Therefore, by the same ODE argument used to show Lemma~\ref{diffeq}, $\hat{\Delta}^n_{r,s}=0$ so $\hG$ is
an $\ainf$ morphism.

Now consider the diagram:
\[\xymatrix{ & \mB & \\ \mA \ar@/^1pc/[ru]^{\sF_1}
\ar@/_1pc/[dr]_{\sF_0}
& \ar@{<-}[l]_-{\iota} \Omega^*(I;\mA) \ar[r]^-{\hG}& \Omega^*(I,\mB), \ar@/_1pc/[lu]_{\text{ev}_1} \ar@/^1pc/[ld]^{\text{ev}_0} \\
& \mB & }\]
where the map $\iota$ maps elements of $\mA$ to constant maps to $\mA$. All
higher order terms of $\iota$ as an $\ainf$ map are 0. By
definition of $\hG$ it is clearly a commutative diagram, so per the discussion
in Section~\ref{sec:diffhtpy}, $\phi= \hG \circ \iota$ and therefore $\hG$ gives
rise to a differential homotopy $\sG^k(\iota^{\tens k})$. 
\end{proof}

\section{Examples} \label{sec:examples}
Let $E^*_{p,q}$ be the cohomological spectral sequence determined by the bigraded complex on
${\Omega^p(M;\Hom(A^{\tens q},A[1-q]))}$ with $d_h=\br{\a^{0,2}}{\cdot}$ and
$d_v=d_{\nabla}$. It follows that $E^2_{p,q}=H^p_{dR}(M;\underline{HH^q(A,A)})$. We
shall compute the cohomology of the total complex in several cases where the $E^2$
term collapses.

Since de Rham cohomology is only defined for globally defined differential
forms on M we shall consider de Rham cohomology with local coefficients as a
version of sheaf cohomology. Using the fact that any two cohomology theories on
$M$ with coefficients in sheaves of \lR-modules over $M$ are uniquely
isomorphic, see for example \cite[p. 184]{warner}, we may use the cohomology
theory that best fits our circumstances. Subsequently, we will
construct a manifold which is homotopy equivalent to $BG$ for any finite group
$G$. This will allow us to transfer the computation of $E^2_{p,q}$ into group
cohomology where it will clearly collapse. Finally, we shall compute the
example of $E^2_{p,q}$ for $M$ homotopy equivalent to a wedge product of
circles, i.e. $BG$ for $G$ a free group with a finite number of generators.

\subsection{Cohomology with local coefficients}
 A theorem of
Eilenberg in \cite{eilenberg1947} tells us that $H^*(X; \E) \approx
H^*_{\text{eq}}(\widetilde{X};E)$ where
$H^*_{\text{eq}}$ are the equivariant cohomology groups and $\widetilde{X}$ is
the universal covering space of $X$ with $\pi_1(X)$ as covering transformations
left operating on E. As noted in \cite{eilenberg1947, brown}, when
$X$ is a $K(G,1)$, $\widetilde{K(G,1)}$ is acyclic and so the augmented
cellular chain complex is a free
resolution of $\R$. 
Equivariant cohomology is that defined on the equivariant cochains
\[C^q_{eq}(\widetilde{BG};E) = \{f \in C^q(\widetilde{BG};E)~|~\delta f(g
\sigma_{q+1})= 
g(\delta f)(\sigma_{q+1})\} \cong \Hom_{G}(C_q(\widetilde{BG}),E)\]
for $g \in G$ and $c_{q+1}$ a $(q+1)$-simplex. Therefore, noting that
$C_*(\widetilde{BG})$ is a chain complex that resolves $\mathbf{R}$, we see that \[C^*(BG; \E) \cong
C^*_{eq}(\widetilde{BG};E) \cong \Hom_G(C_*(\widetilde{BG}),E) = C^*(G;E),\] and thus
\[H^*(BG;\E) \cong H^*(G,E).\]

\subsection{Making $BG$ a manifold} \label{bgasmfld}

Let $\Gamma$ be a finite group. First, recall that every finite group
$\Gamma$ of order $k$ embeds in the unitary group $U(k)$. This follows
by noting that every finite group has a faithful representation in
$GL(k)$ given by the permutation representation and that every finite
subgroup of $GL(k)$ is conjugate to a subgroup of $U(k)$ (see, for
example \cite[\S 9.2]{artin}). Now $U(k)$ is a compact Lie group and
certainly a manifold. The Grassmannian $G(k,n+k) = U(n+k)/(U(n)\times
U(k))$ is a smooth compact manifold and $BU(k)= \lim_{n\to \infty}G(k,n+k)$. The group $U(k)$ acts
freely on $U(n+k)$ and on $U(n+k)/U(n)$, so $\Gamma \subset U(k)$ also
acts freely on both spaces. Consider the
space $U(n+k)/U(n)$ as the space of orthonormal families of $k$
vectors in $\mathbf{C}^{n+k}$. Thus we have a fibre bundle
\[\xymatrix{U(n-1+k)/U(n-1) \ar[r] & U(n+k)/U(n) \ar[d] \\
   & S^{2(n+k)-1}.}\]
Therefore, again taking $n \to \infty$, we see that $EU(k)= \lim_{n\to \infty} U(n+k)/U(n)$ is contractible and has a free $\Gamma$ action so $EU(k)/\Gamma$ is a classifying space for $\Gamma$ (see, for example \cite{hatcher}).

For brevity, let us denote submanifold $U(n+k)/U(n)\times \Gamma \subset G(k,n+k)$ by $M^n$. Therefore, we have constructed a series of
manifolds \[M^1 \hookrightarrow M^2
\hookrightarrow M^3 \hookrightarrow \cdots \hookrightarrow M^k \hookrightarrow \cdots\] using the natural inclusions with the property that the
limit space $EU(k)/\Gamma=\bigcup_{i=1}^{\infty}M^i$ has the same homotopy type as $BG$.

We want to be able to talk
about $\Omega^*(B\Gamma)$, and thus need a notion of a path. Define ``a
path in $B\Gamma$'' as a path in $M^i$ for some $i$. Then
\[\Omega^*(B\Gamma) = \{ \theta_i \in \Omega^*(M^i)_{i\geq 0}~|~
\theta_i|_{M^{i-1}} = \theta_{i-1}\}.\] This is an inverse limit system which
trivially satisfies the Mittag-Leffler condition
because $M^{i-1} \subset M^i$ \cite[p. 191]{hartshorne}. We can define a differential and a wedge product on these differential
forms. Given a differential form $\theta = \lim\limits_{\leftarrow} \theta_i \in \Omega^*(B\Gamma)$, define
\[d\theta = \lim\limits_{\leftarrow} d \theta_i \hspace{.2in} \in
\Omega^*(B\Gamma).\] Similarly, given two differential forms $\omega =
\lim\limits_{\leftarrow} \omega_i$ and $\eta = \lim\limits_{\leftarrow} \eta_i$
in $\Omega^*(B\Gamma)$, define the wedge product \[\omega \wedge \eta =
\lim\limits_{\leftarrow} (\omega_i \wedge \eta_i) \hspace{.2in}\in
\Omega^*(B\Gamma).\] With these maps, we can consider the de Rham
cohomology of $B\Gamma$.

Given a suitable topology on $B\Gamma$, we know
that \[H^*_{dR}(B\Gamma; \E) \cong H^*_{\Delta}(B\Gamma; \E).\] Thus we can
compute cohomologies using the
simplicial cohomology of $B\Gamma$ with twisted coefficients in $HH^*(A,A)$ where we consider $B\Gamma$ as the
standard simplicial complex generated from the universal cover
$E\Gamma$. In this case, there is one vertex, *, in $B\Gamma$. The
simplices of $B\Gamma$ can be described using the bar notation
$[g_1|g_2|\cdots |g_n]$. In this notation, the boundary simplices of
$[g_1|g_2|\cdots|g_n]$ are $[g_2|\cdots|g_n],~[g_1|\cdots|g_{n-1}],$
and $[g_1|\cdots|g_ig_{i+1}|\cdots|g_n]$ for $i=1,\ldots,n-1$. Since
this is independent of topology, we shall not specify one.

\subsection{Finite groups} \label{ex:finitegroups}
Let $\Gamma$ be a finite group. Construct the manifold $B\Gamma$ as in
Section~\ref{bgasmfld}. Let \underline{A} 
be an $A$ local system on $B\Gamma$ where $A$ is a real vector space that
is an associative algebra under the multiplication map
$\mu^2=\a^{0,2}$. Consider the total complex
$\Omega^*(B\Gamma ;\Hom(A^{\tens *},A[1-*]))$ with differential
$d_{\nabla}+[\a^{0,2},\cdot]$. As a spectral sequence with
$E^0_{m,n}=\Omega^m(B\Gamma ;\Hom(A^{\tens n},A[1-n]))$, we have
\[E^2_{m,n}=H^m_{dR}(B\Gamma ;\underline{HH^n(A,A)})=H^m_{\text{simp}}(B\Gamma ;\underline{HH^n(A,A)}).\]  By the discussion above, it is
sufficient to calculate $H^*(\Gamma ;HH^*(A,A))$. However,
by \cite[pg. 117]{maclane}, because $HH^n(A,A)$ is
always a real vector space and therefore a divisible abelian group with no
elements of finite order \[H^p(\Gamma ;HH^{n}(A,A)) = \begin{cases} HH^n(A,A)^{\Gamma}
& \text{if }p=0 \\
0 &\text{if } p\neq 0.\end{cases}\]
Thus, $E^2$ degenerates to a single nonzero column so $E^2 \simeq E^{\infty}$.

\subsection{Finitely generated free group}
\label{sec:nonabex}
Let $G$ be a finitely generated nonabelian free group with $r$ generators, let $S$ be a set
with $r$ elements and let $Y=\vee_{s \in S} S^1_s$. Then $Y$ is clearly a
$K(G,1)$ because 
\[\pi_1(Y) = G \] and \[\pi_k(Y) = \bigoplus_{s\in
S}\pi_k(S^1_s) = 0 \text{ for }k \geq 2.\]
It is also straightforward to see that we may consider $Y$
as a manifold with only one inflationary step of the type used in
Section~\ref{bgasmfld}. Therefore we again
have \[H^p(Y; \underline{HH^q(A,A)}) \cong H^p(G; HH^q(A,A))\] which we can
compute using the resolution:
\begin{equation}
\xymatrix{0 \ar[r] & \R G^{(S)} \ar[r]^-{\partial}& \R G \ar[r]^-{\epsilon}
&\R \ar[r] & 0}  
\end{equation}
where $\R G^{(S)}$ has basis $t_s$ corresponding to the oriented 1-simplex
mapping to $S^1_s$, $\R G$ has basis $x$ corresponding to the base point and
$\partial(t_s) = (1-g_s)x$ because we must translate the endpoints of
$\Delta^1$ to the same point before summing.

Now the cohomology $H^*(G;HH^q(A,A))$ is the cohomology of the complex
\[\xymatrix{HH^q(A,A) \ar[r]^-{\delta} &  \bigoplus_{s \in S}HH^q(A,A)_s \ar[r]
& 0 \ar[r] & \cdots,}\] with $(\delta u) = \oplus_{s\in S}(1-g_s)u$ . Therefore,
\[ H^p(G;HH^q(A,A)) = \begin{cases} HH^q(A,A)^G & \text{ if } p=0 \\
 \text{Coker}(\delta) &\text{ if }p=1 \\
 0 & \text{ otherwise.}\end{cases}\]
Hence, $E^2_{p,q}$ has two nonzero columns, but that still means that
 $E^2_{p,q}\cong E^{\infty}_{p,q}$.  

\section{Transferring Maurer-Cartan solutions between $CC^*(A,A)^{\Gamma}$ and
 $\Omega^*(B\Gamma;\g)$}   

 \begin{proposition} \label{prop:GinGout} Let $\Gamma$ be a finite group. Then
 \[HH^*(A,A)^{\Gamma} \cong H^*(\g_x^{\Gamma}).\]
 \end{proposition}
 \begin{proof}
  First, note that $H^*(\g_x^{\Gamma})$ is
 isomorphic to $H^*(CC^*(A,A)^{\Gamma})$ so we shall use them interchangeably.
 It is clear that $H^q(\g_x^{\Gamma}) \subset HH^q(A,A)^{\Gamma}$, so it only
 remains to check the opposite inclusion.
Let $[f] \in HH^q(A,A)^{\Gamma}$ for $f \in CC^q(A,A)$ closed. Consider the
 element \[\overline{f} = \frac{1}{|\Gamma|} \sum_{g \in \Gamma} gf.\] Now,
 $\overline{f}$ is a closed element of $CC^q(A,A)^{\Gamma}$ by construction and because $[f]$ represents
 a $\Gamma$-invariant class it is clear that $[\overline{f}]=[f] \in
 HH^q(A,A)^{\Gamma}$.
Now consider changing $f$ by a coboundary $d\omega$. The averaging process then
 shows that \[\overline{f+d\omega} = \frac{1}{|\Gamma|} \sum_{g\in \Gamma} gf +
 gd\omega = \frac{1}{|\Gamma|} \sum_{g\in \Gamma} gf +
 dg\omega\] since $g\mu^2 = \mu^2(g\tens g)$ and therefore $[\overline{f}] = [\overline{f}+d\overline{\omega}]$
 so we can also choose the coboundary representatives in $CC^{q-1}(A,A)^{\Gamma}$. \end{proof}

\begin{lemma}\label{lem:cctoomega} Let $M$ be a manifold with a basepoint. Consider $CC^*(A,A)^{\pi_1(M)}$ as a
 differential graded Lie algebra under the Hochschild differential and the Gerstenhaber bracket.
 Then there
 is a dg Lie map $\eta:CC^*(A,A)^{\pi_1(M)} \to \Omega^*(M;\g)$ defined by $\eta(f) \mapsto$ \{the
 constant 0-form with value $f$ at the basepoint\} and extending this to all of $\g$ by
 parallel transport.  \end{lemma} \begin{proof} 
 There are no difficulties with nontrivial loops because
 $f \in CC^*(A,A)^{\pi(M)}$. Also, $\eta$ clearly commutes with the differentials
 because the DeRham differential on $\Omega^*(M;\g)$ is zero on constant forms
 leaving only the Hochschild differential in each case and $\eta$ commutes with
 the bracket because the wedge of constant 0-forms is a constant 0-form and
 will not change any signs because it is of even degree.
 \end{proof}

 \begin{corollary}
 When $M= B\Gamma$ for a finite group $\Gamma$, then $CC^*(A,A)^{\Gamma}$ is
 quasi-isomorphic to $\Omega^*(M;\g)$.
 \end{corollary}
 \begin{proof} The map $\eta$ constructed in Lemma~\ref{lem:cctoomega} induces
 a map \[\overline{\eta}: H^*(CC^*(A,A)^{\Gamma}) \to H^*(\Omega^*(M;\g)),\] but by
 Section~\ref{ex:finitegroups} we know that $H^*(\Omega^*(M;\g))\cong
 HH^*(A,A)^{\Gamma}$ and by Proposition~\ref{prop:GinGout} we know that
 $H^*(CC^*(A,A)^{\Gamma}) = HH^*(A,A)^{\Gamma}$. \end{proof}

 Let $\h \subset CC^*(A,A)^{\Gamma}$ be the sub-dg Lie algebra of functions
 with negative internal degree in $A$, i.e.
 \[\h = \{f\in CC^*(A,A)^{\Gamma}~|~f \in \Hom(V^{\tens *},V[-n]), n \geq 1\}.\]
 There is a natural decreasing filtration $L_k\h$ for $k \geq 1$ given by \[L_k\h = \{f\in
 CC^*(A,A)^{\Gamma}~|~f \in \Hom(V^{\tens *},V[-n]), n \geq k\}.\] Thus, $L_1\h=\h$ and $\mu^2 \in \Hom(V\tens V, V[0])$ so for
 $f \in \Hom(V^{\tens \ell},V[-m])$ we
 have \[\br{\mu^2}{f} \in \Hom(V^{\ell+1},V[-m])\] and thus
 $d(L_m\h) \subset L_m\h$. Lastly, $\br{\h^m}{\h^n} \subset \h^{m+n}$ by the
 additivity of degrees so if we include a formal degree 0 parameter $\hbar$
 so that $F_k\h =
 L_i\h\hbar^k$, then $\h$ with filtration $F_{\bullet}$ is a filtered pronilpotent dg Lie algebra.

 Let $\j \subset \Omega^*(M;\g)$ be, similarly, the sub-dg Lie algebra of
 functions with negative internal degree strictly less than $-1$. Our checks in
 Proposition~\ref{prop:Oisdgla} show that \[L_k\j
 = \{\Omega^*(M;\Hom^{\bullet}(\underline{V}^{\tens \bullet},\underline{V}[-k])\}, \hspace{.2in}k \geq
 1\]  and $F_k\j = L_k\j\hbar^k$ make $\j$ a filtered pronilpotent dg Lie algebra.

\begin{definition}For $\mathfrak{f}$ a dg Lie algebra, let $MC(\mathfrak{f})$ denote the set of solutions to the Maurer-Cartan
 equation in $\mathfrak{f}$. An element $\a \in MC(\mathfrak{f})$ will be
 called a Maurer-Cartan element.\end{definition} 

\begin{lemma} Let $\h$ and $\g$ be filtered pronilpotent dg Lie
 algebras. Suppose that $\Phi:\h \to \g$ is a 
 filtered quasi-isomorphism between them which
  means that $\Phi$ induces quasi-isomorphisms of chain complexes
 $F_r\h/F_{r+1}\h \to F_r\g/F_{r+1}\g$ for any $r$. Then $\Phi$ induces a
 bijection between equivalence classes of Maurer-Cartan elements.
 \end{lemma}
 \begin{proof}
 This is a special case of the isomorphism between deformation functors of
 Section~4.4 of \cite{defquant}. A similar formulation in terms of filtered
 dglas can be found in Lemma~2.2 of \cite{hmsg2c}. \end{proof}

\begin{corollary} \label{cor:mchtomco} The map 
$\eta$ defined in Lemma~\ref{lem:cctoomega} is a filtered quasi-isomorphism
between $\h \subset CC^*(A,A)^{\Gamma}$ and $\j \subset \Omega^*(M;\g)$. Thus
the map $\a \mapsto \eta(\a)$ induces a bijection between equivalence classes
of Maurer-Cartan elements.\end{corollary}

Let $\a \in MC(\h)$ be a Maurer-Cartan element. Then, since
$\a \in CC^1(A,A)^{\Gamma} \cap \h$ we see that $\a= \sum_{k=1}^{\infty}\a^{k+2}\hbar^k$
where $\a^k \in \Hom(V^{\tens k},V[2-k])$. Likewise, we have Maurer-Cartan
elements $\ha \in \j$. The total degree one constraint in $\Omega^*(M;\g)$
requires that $\ha^{m,n} = 0$ for $m+n \leq 2$ if a solution is to be in
$\j$. Therefore, $\ha^{0,0}
=0, \ha^{0,1}=0, \ha^{0,2}=0, \ha^{1,0}=0, \ha^{1,1}=0$, and $\ha^{2,0}=0$. In
particular, Maurer-Cartan elements in $\j$ will satisfy all the assumptions
that we applied in Section~\ref{sec:families}.

\subsection{Homotopies of Maurer-Cartan elements}
There is a natural Lie algebra homomorphism from $\j^0$ to the space of affine
vector fields on $\j^1$ which associates to $\ga \in \j^0$ the infinitesimal
gauge transformation \[\a \mapsto -d_{\j}\ga+\br{\ga}{\a}.\] Using the
Baker-Campbell-Hausdorff formula, and the fact that we are working in
pronilpotent dg Lie algebras there is a group action on the set of
Maurer-Cartan elements by $\exp(\j^0)$. For $\ga \in \j^0$ and $\a\in MC(\j)$, let us denote the action of $\exp(\ga)$
on $\a$ by $\circledast$ while the infinitesimal action of $\ga$ on $\a$ is
simply denoted by $\circledcirc$.

Let $\a_0$ and $\a_1$ be two equivalent Maurer-Cartan elements in $\j^1$. We shall
construct a homotopy between them.
Since $\a_0$ is equivalent to $\a_1$, there exists a $\ga \in \j^0$ so that
\[\a_1 = \exp(\ga)\circledast\a_0= \a_0 -d_{\j}\ga +
\br{\ga}{\a_0}-\frac{1}{2!}\br{\ga}{d_{\j}\ga}+\frac{1}{2!}\br{\ga}{\br{\ga}{\a_0}} +\cdots.\] For any element
$\ga \in \j^0$, it follows that $t\ga \in \j^0$ for $t \in [0,1]$. Consider the
one parameter family of Maurer-Cartan
elements \[\a_t=\exp(t\ga)\circledast\a_0.\] This family has the property that
the derivative at each point $t \in [0,1]$ is defined in terms of $\ga$ and the
value at the point $\at$.
\begin{align*}  \der{\at}{t} &= \ga\circledcirc\left(\exp(t\ga)\circledast
  \a_0\right) \\
  &= -d_{\j} \ga + \br{\ga}{\at}
  \end{align*}
  
Consider the space $\Omega^*(I_t;\j)$ as a dg Lie algebra with differential
$D=\der{}{t}dt+d_{\j}$ and the bracket induced by the bracket on
$\j$ and the wedge product of forms. Then $\a_t+\ga dt$ is a
Maurer-Cartan element because
\begin{align*}
  D(\at+\ga dt) + \frac{1}{2}\br{\at+\ga dt}{\at+\ga dt}& = d_{\j} (\at) +
  \frac{1}{2}\br{\at}{\at} \\ & +\left(\der{\at}{t} + d_{\j}(\ga) +
  \br{\at}{\ga}\right)dt.
  \end{align*}
As the bracket is induced, the first line is 0 because we know that $\at$ is a
Maurer-Cartan element in $\j$ for all $t$. Likewise, the second line is zero
because \[
\der{\at}{t}+d_{\j}\ga+\br{\at}{\ga}=-d_{\j}\ga+\br{\ga}{\at}+d_{\j}\ga-(-1)^{||\ga||||\at||}\br{\ga}{\at}
= 0.\]

\subsection{Map from $\Omega^*(I;\j) \to \Omega^{*}(I\times M;\g^{\text{neg}})$}
Define $\g^{\text{neg}}\subset \g$ to be all elements of $\g$ with negative
internal degree. Thus, we can recall that $\Omega^*(M;\g^{\text{neg}})=\j$. We wish to say that the map \begin{align*}\iota:\Omega^*(I;\j) &\to \Omega^*(I\times
M;\g^{\text{neg}})\\ f+gdt &\mapsto f+g\wedge dt\end{align*} will allow us to transfer
Maurer-Cartan elements to homotopies of Maurer-Cartan elements. There is clearly a map
$\pi:\Omega^*(I\times M;\g^{\text{neg}}) \to \Omega^*(I;\j)$ given by restricting forms
on $I\times M$ with values in $\g^{\text{neg}}$ to forms on $I$ with values in $\j$.

First, observe that $\ga \in \j^0$ is independent of $t$. Second, note that
$\at=\exp(t\ga)\circledast \a_0$, and therefore the degree in $t$ rises
with each included term of $\ga$.  We want $\ga \in \j^0$, so that means $\ga \in \Omega^m(M;\g^{-m})$ or in other
words, $\ga^{0,0}=\ga^{0,1}=\ga^{1,0}=0$. Thus,
\begin{align}\label{eqn:gamma}\nonumber \ga &= (\ga^{0,2}+\ga^{1,1}+\ga^{2,0})\hbar +\\&\hspace{.23in}
(\ga^{0,3}+\ga^{1,2}+\ga^{2,1}+\ga^{3,0})\hbar^2 +\\
\nonumber  &\hspace{.23in}
  (\ga^{0,4}+\ga^{1,3}+\ga^{2,2}+\ga^{3,1}+\ga^{4,0})\hbar^3 +
  \cdots
  \end{align}
Let us filter $\at=\exp(t\ga)\circledast \a_0$ by our filtration $F_k\j$. The first
few terms are:
\begin{align*}\at
  &= -d_{\j}(t\ga)+\a_0 & \hspace{-.2in}\mod F_2\j \\ \at&=
  -d_{\j}(t\ga)+\a_0+t\br{\ga}{\a_0}-\frac{t^2}{2!}\br{\ga}{d_{\j}\ga} &\hspace{-.2in}\mod F_3\j \\
 & \hspace{.4in}\vdots \hspace{1in}\vdots \hspace{1in}\vdots \\ \at &=
  -d_{\j}(t\ga)+\a_0+t\br{\ga}{\a_0} + \cdots
  - \frac{t^k}{k!}[\underbrace{\ga,[\ga,[\cdots,[\ga}_{k-1},d_{\j}\ga]\cdots]]]
  &\hspace{-.2in}\mod F_{k+1}\j \end{align*} where each term contributes only a finite number
  of components because of the filtration as shown in (\ref{eqn:gamma}). Thus,
  the pronilpotence of $\j$ takes care of the convergence of the map.

  \begin{theorem} \label{thm:alltrivial} For a finite group $\Gamma$,
  every homotopy $\Gamma$ action on an $\ainf$ algebra $\mA$ has class
  representatives $\sF_g:\mA\to \mA$ for all $g \in \Gamma$ which comprise a
  strict action. Therefore, $\sF_g \circ \sF_h = \sF_{gh}$ and $\sF_e=\id$.\end{theorem}
\begin{proof}  
  Let $p$ be a closed loop in $M\cong K(\Gamma,1)$ based at a point $x$ where
  $[p]=g \in \Gamma$. The homotopy group action of $\Gamma$ on $\mA_x$ is defined by the
  actions of the generators $g$ on $\mA_x$. In Figure~\ref{mchtpyclasses} we
  see that the loop $p$ defines a cylinder in $I\times M$.
\begin{center}
\begin{figure}[htbp]
\includegraphics[width=5in]{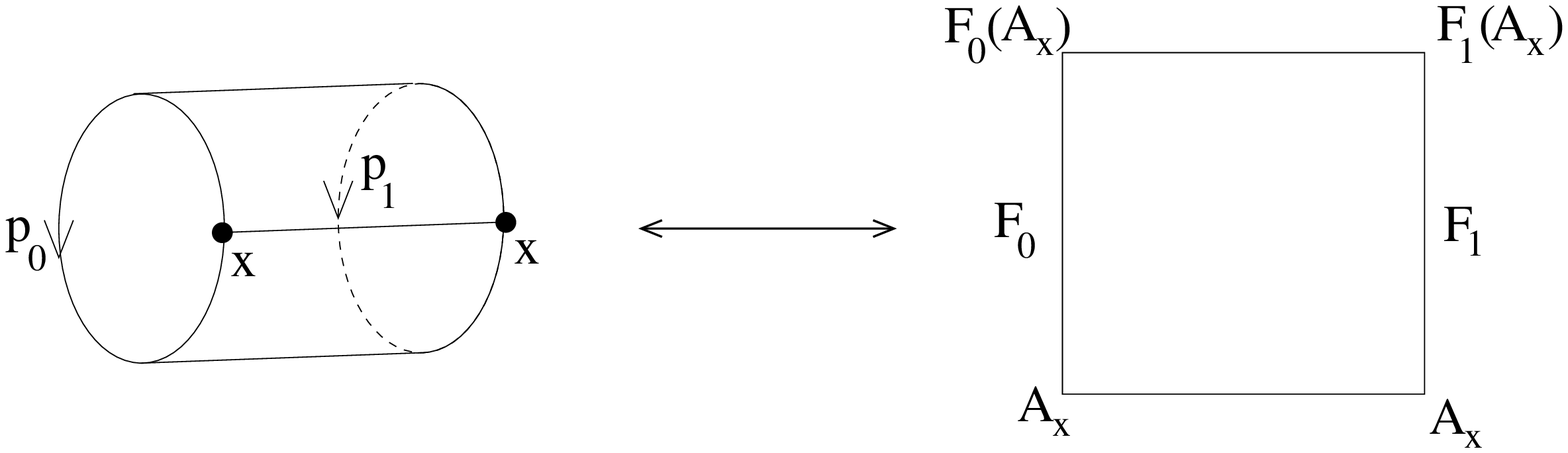}
\caption{The correspondence between the cylinder $p\times I$ and the homotopy
  between $\sF_0$ and $\sF_1$.\label{mchtpyclasses}}
\end{figure}
\end{center}
Let $[\sF_g]$ denote the class of $\ainf$ endomorphisms of $\mA_x$ that
  correspond to the $\Gamma$ homotopy group action. The $\ainf$ morphism
  $\sF_{p,0} \in [\sF_g]$ is defined by integrating $\a_0\in MC(\j)$ according to
  Theorem~\ref{thm:sFn} while $\sF_{p,1} \in [\sF_g]$ is defined by integrating
  $\a_1 \in MC(\j)$ accordingly. Integrating over the square $I\times I_s$
  where $I_s$ corresponds to traversing $p$ gives a homotopy
  $T:\sF_{p,0} \to \sF_{p,1}$ by applying Theorem~\ref{thm:hGn}. Consider the the
  $\ainf$ endomorphism $\sF_g = \sG_{0\to 1} \circ \sF_{p,1} \circ \sG_{1\to 0}$
  where $\sG_{a\to b}$ consists of integrating the appropriate terms of
  $\a_t+\ga dt \in MC(\Omega^*(I\times M;\g))$ along the path $\{x\}\times
  I$. Now consider a second path $q:I_s \to M$ with $[q]=h \in \Gamma$ and
  construct $\sF_h = \sG_{0\to 1}\circ \sF_{q,1}\circ \sG_{1,0}$ in the same
  manner. By construction, it is clear that $\sF_g \in [\sF_g]$ and $\sF_h \in
  [\sF_h]$.  Let \[pq = \begin{cases} p(2s) & 0\leq s \leq \frac{1}{2} \\
  q(2s-1) & \frac{1}{2}\leq s \leq 1 \end{cases}.\] Now by applying
  Corollary~\ref{cor:strictinverses} to eliminate $\sG_{1\to 0}\circ \sG_{0\to
  1}$ and Lemma~\ref{compositionpd} to combine $\sF_{p,1}\circ \sF_{q,1}$ we have\begin{align*}\sF_g\circ \sF_h &= \sG_{0\to
  1}\circ \sF_{p,1}\circ \sG_{1\to 0}\circ \sG_{0\to
  1}\circ \sF_{q,1}\circ \sG_{1,0} \\
  &= \sG_{0\to
  1}\circ \sF_{p,1}\circ \sF_{q,1}\circ \sG_{1,0} \\
   &= \sG_{0\to
  1}\circ \sF_{pq,1}\circ \sG_{1,0} \\
&= \sF_{gh} \in [\sF_{gh}],\end{align*}
as desired.

\end{proof}

\end{document}